\documentclass{amsart}
\usepackage{amsfonts, amsmath, latexsym, epsfig}
\usepackage[norelsize]{algorithm2e}
\usepackage{amssymb}
\usepackage{color}
\usepackage{epsf}
\usepackage{url}
\usepackage{tikz}
\tikzset{dynnode/.style={circle, draw, thick, inner sep=1pt, minimum size=0.8cm, font=\small}}

\newcommand{\RR}{\ensuremath{\mathbb{R}}}

\newcommand{\ZZ}{\ensuremath{\mathbb{Z}}}

\newtheorem{proposition}{Proposition}
\newtheorem{theorem}{Theorem}

\newtheorem{lemma}{Lemma}

\newtheorem{remark}{Remark}
\newtheorem{definition}{Definition}
\def\QuotS#1#2{\leavevmode\kern-.0em\raise.2ex\hbox{$#1$}\kern-.1em/\kern-.1em\lower.25ex\hbox{$#2$}}

\newcommand{\MC}{{\mathcal C}}

\DeclareMathOperator{\Stab}{Stab}

\begin{document}

\author{Mathieu Dutour Sikiri\'c}
\address{Mathieu Dutour Sikiri\'c, MSM programming, Karlovacka cesta 28B, 10452 Klinca selo}
\email{mathieu.dutour@gmail.com}

\title[Minimal norm of regular vectors]{Minimal norm of regular vectors in the root lattice of a crystallographic root system}

\date{}

\begin{abstract}
For a crystallographic root system $\Phi$ with Weyl group $W$ and root lattice $L = \ZZ\Phi$,
we determine the minimum norm of a vector of $L$ having trivial stabilizer in $W$,
and we prove that the vectors attaining it form a single orbit of $W$.
Although this is a minimization of a quadratic function over the lattice points of a cone,
the minimizer in the dominant chamber is always unique:
it is $\rho + \lambda$, where $\rho$ is the Weyl vector and $\lambda$ is $0$ or the minuscule weight
in the class of $\rho$ modulo $L$.
In particular, the optimal vector is the Weyl vector exactly when the Weyl vector lies in $L$.
We give explicit formulas for the families $\mathsf{A}_n$, $\mathsf{B}_n$, $\mathsf{C}_n$ and $\mathsf{D}_n$, with elementary proofs,
and the values for the exceptional types $\mathsf{E}_6$, $\mathsf{E}_7$, $\mathsf{E}_8$, $\mathsf{F}_4$ and $\mathsf{G}_2$.
Finally, we explain the coincidences between the minimal norms by the folding of Dynkin diagrams.
\end{abstract}

\maketitle

\section{Introduction}

Let $\Phi$ be a crystallographic root system in a Euclidean space $V$
with Weyl group $W = W(\Phi)$, and let $L = \ZZ\Phi$ be its root lattice.
We refer to Humphreys~\cite[Chapters~1 and~2]{Humphreys} for the definitions and basic properties of
root systems, Weyl groups and root lattices;
crystallographic root systems, root lattices and weight lattices are treated in~\cite[\S 2.9]{Humphreys}.
For the simply laced types, with roots of norm $2$, the root lattices are
$\mathsf{A}_n$, $\mathsf{D}_n$, $\mathsf{E}_6$, $\mathsf{E}_7$ and $\mathsf{E}_8$.
For the types $\mathsf{B}_n$ and $\mathsf{C}_n$ we use the standard realizations in $\RR^n$:
\begin{equation*}
\Phi(\mathsf{B}_n) = \{\pm e_i, \pm e_i \pm e_j\}
\quad\text{and}\quad
\Phi(\mathsf{C}_n) = \{\pm 2 e_i, \pm e_i \pm e_j\},
\end{equation*}
whose root lattices are $\ZZ^n$ and $\mathsf{D}_n$ respectively.
In both cases the Weyl group is the group of all signed permutations of the coordinates,
which for $\mathsf{C}_n$ is larger than $W(\mathsf{D}_n)$.
For $\mathsf{G}_2$ we take, in the plane $x_1 + x_2 + x_3 = 0$ of $\RR^3$,
the short roots $e_i - e_j$ and the long roots $\pm(2e_i - e_j - e_k)$ with $\{i,j,k\} = \{1,2,3\}$;
the root lattice is $\mathsf{A}_2$ and $W(\mathsf{G}_2)$ is the dihedral group of order $12$.
For $\mathsf{F}_4$ we take, in $\RR^4$, the short roots $\pm e_i \pm e_j$ and the long roots
$\pm 2e_i$ and $(\pm 1, \pm 1, \pm 1, \pm 1)$;
the root lattice is $\mathsf{D}_4$ and $W(\mathsf{F}_4)$ has order $1152$.
The root lattice alone therefore does not determine the problem considered below:
$\mathsf{C}_n$ and $\mathsf{D}_n$ share the lattice $\mathsf{D}_n$, and $\mathsf{F}_4$ and $\mathsf{D}_4$ share the lattice $\mathsf{D}_4$,
but with different Weyl groups.

The reflection hyperplanes
$H_\alpha = \{x \in V : \langle x, \alpha \rangle = 0\}$, for $\alpha \in \Phi$,
partition $V$ into open \emph{Weyl chambers}.
Since the stabilizer of a vector $v \in V$ is generated by the reflections it contains~\cite[Theorem~1.12]{Humphreys},
a lattice vector $v \in L$ has trivial stabilizer in $W$ if and only if
$v$ lies in the interior of some Weyl chamber, equivalently,
$\langle v, \alpha \rangle \neq 0$ for all roots $\alpha \in \Phi$.
Such vectors are called \emph{regular}.

We consider the following optimization problem.

\begin{definition}
For a root system $\Phi$, define
\begin{equation*}
\operatorname{Triv}(\Phi) = \min \left\{ \lVert v \rVert^2 : v \in \ZZ\Phi, \; \Stab_{W(\Phi)}(v) = \{1\} \right\}.
\end{equation*}
\end{definition}

By symmetry, we may restrict the search to a single closed Weyl chamber.
Fixing a system of simple roots $\alpha_1, \ldots, \alpha_n$, the dominant chamber is
$\MC = \{ x : \langle x, \alpha_i \rangle \geq 0 \text{ for all } i \}$,
and we seek a lattice point in its interior with minimum norm.
In the realizations above the inner product takes integral values on $L$, so that this is
\begin{equation}
\label{eq:iqp}
\operatorname{Triv}(\Phi) = \min \left\{ v G v^T : v \in \ZZ^n, \; v G \alpha_i^T \geq 1 \text{ for all } i \right\},
\end{equation}
where $G$ is the Gram matrix of the simple roots and $v$ is expressed in the basis of simple roots.

This is an integer quadratic programming problem that can be solved computationally,
and a priori nothing guarantees that its optimal solution is unique.
Our main result, Theorem~\ref{thm:minuscule}, shows that it always is:
the optimal vector in the dominant chamber is $\rho + \lambda$,
where $\rho$ is the Weyl vector and $\lambda$ is the minimal dominant weight in the class of $\rho$ modulo $L$,
which is $0$ or a minuscule weight.
For the families $\mathsf{A}_n$, $\mathsf{B}_n$, $\mathsf{C}_n$ and $\mathsf{D}_n$, we also give elementary proofs of closed-form solutions
(Theorems~\ref{thm:An} and~\ref{thm:BCD}).
For the exceptional types, the solutions were first obtained computationally
and are reported in Section~\ref{sec:exceptional}.
In Section~\ref{sec:weyl} we compare the optimal vectors with the Weyl vector and prove Theorem~\ref{thm:minuscule},
and in Section~\ref{sec:folding} we explain the equalities between the minimal norms
of different root systems by folding.
The results are summarized in Table~\ref{tab:results}.

\begin{table}[h]
\centering
\renewcommand{\arraystretch}{1.3}
\begin{tabular}{llll}
$\Phi$ & $\operatorname{Triv}(\Phi)$ & $\lambda$ & Reference\\
\hline
$\mathsf{A}_n$, $n$ even & $\frac{k(k+1)(2k+1)}{3}$, $k = \frac{n}{2}$ & $0$ & Theorem~\ref{thm:An}\\
$\mathsf{A}_n$, $n$ odd & $\frac{k(k+1)(2k+1)}{3}$, $k = \frac{n+1}{2}$ & $\omega_k$ & Theorem~\ref{thm:An}\\
$\mathsf{B}_n$ & $\frac{n(n+1)(2n+1)}{6}$ & $\omega_n$ & Theorem~\ref{thm:BCD}\\
$\mathsf{C}_n$, $n \equiv 0, 3 \pmod{4}$ & $\frac{n(n+1)(2n+1)}{6}$ & $0$ & Theorem~\ref{thm:BCD}\\
$\mathsf{C}_n$, $n \equiv 1, 2 \pmod{4}$ & $\frac{n(n+1)(2n+1)}{6} + 2n + 1$ & $\omega_1$ & Theorem~\ref{thm:BCD}\\
$\mathsf{D}_n$, $n \equiv 0, 1 \pmod{4}$ & $\frac{(n-1)n(2n-1)}{6}$ & $0$ & Theorem~\ref{thm:BCD}\\
$\mathsf{D}_n$, $n \equiv 2, 3 \pmod{4}$ & $\frac{(n-1)n(2n-1)}{6} + 2n - 1$ & $\omega_1$ & Theorem~\ref{thm:BCD}\\
$\mathsf{E}_6$ & $78$ & $0$ & Figure~\ref{fig:exceptional}\\
$\mathsf{E}_7$ & $228$ & $\omega_7$ & Figure~\ref{fig:exceptional}\\
$\mathsf{E}_8$ & $620$ & $0$ & Figure~\ref{fig:exceptional}\\
$\mathsf{F}_4$ & $78$ & $0$ & Figure~\ref{fig:exceptional}\\
$\mathsf{G}_2$ & $14$ & $0$ & Figure~\ref{fig:exceptional}
\end{tabular}
\caption{The minimal norm $\operatorname{Triv}(\Phi)$ for each irreducible root system $\Phi$,
the weight $\lambda$ such that the optimal vector in the dominant chamber is $\rho + \lambda$ (Theorem~\ref{thm:minuscule}),
and where the optimal vector is described.
The fundamental weights are numbered as in Bourbaki;
$\lambda = 0$ means that the optimal vector is the Weyl vector $\rho$.}
\label{tab:results}
\end{table}

\begin{remark}
\label{rem:reducible}
The problem reduces to irreducible root systems.
Let $\Phi = \Phi_1 \cup \cdots \cup \Phi_r$ be the decomposition of $\Phi$ into irreducible components,
which are pairwise orthogonal.
Then $L = L_1 \oplus \cdots \oplus L_r$ with $L_j = \ZZ\Phi_j$, and $W = W(\Phi_1) \times \cdots \times W(\Phi_r)$,
each factor acting on the span of its own component and trivially on the others.
For $v = v_1 + \cdots + v_r$ with $v_j \in L_j$, we have
$\Stab_W(v) = \Stab_{W(\Phi_1)}(v_1) \times \cdots \times \Stab_{W(\Phi_r)}(v_r)$
and $\lVert v \rVert^2 = \sum_j \lVert v_j \rVert^2$.
Hence
\begin{equation*}
\operatorname{Triv}(\Phi) = \operatorname{Triv}(\Phi_1) + \cdots + \operatorname{Triv}(\Phi_r),
\end{equation*}
and the optimal vectors of $\Phi$ are the sums of optimal vectors of the components.
In particular, the optimal solution is again unique up to the action of $W$,
and the optimal vector in the dominant chamber is $\rho + \lambda$,
where $\rho$ and $\lambda$ are the sums of the corresponding vectors of the components.
In the rest of the paper we assume that $\Phi$ is irreducible.
\end{remark}

\section{Exceptional root systems}
\label{sec:exceptional}

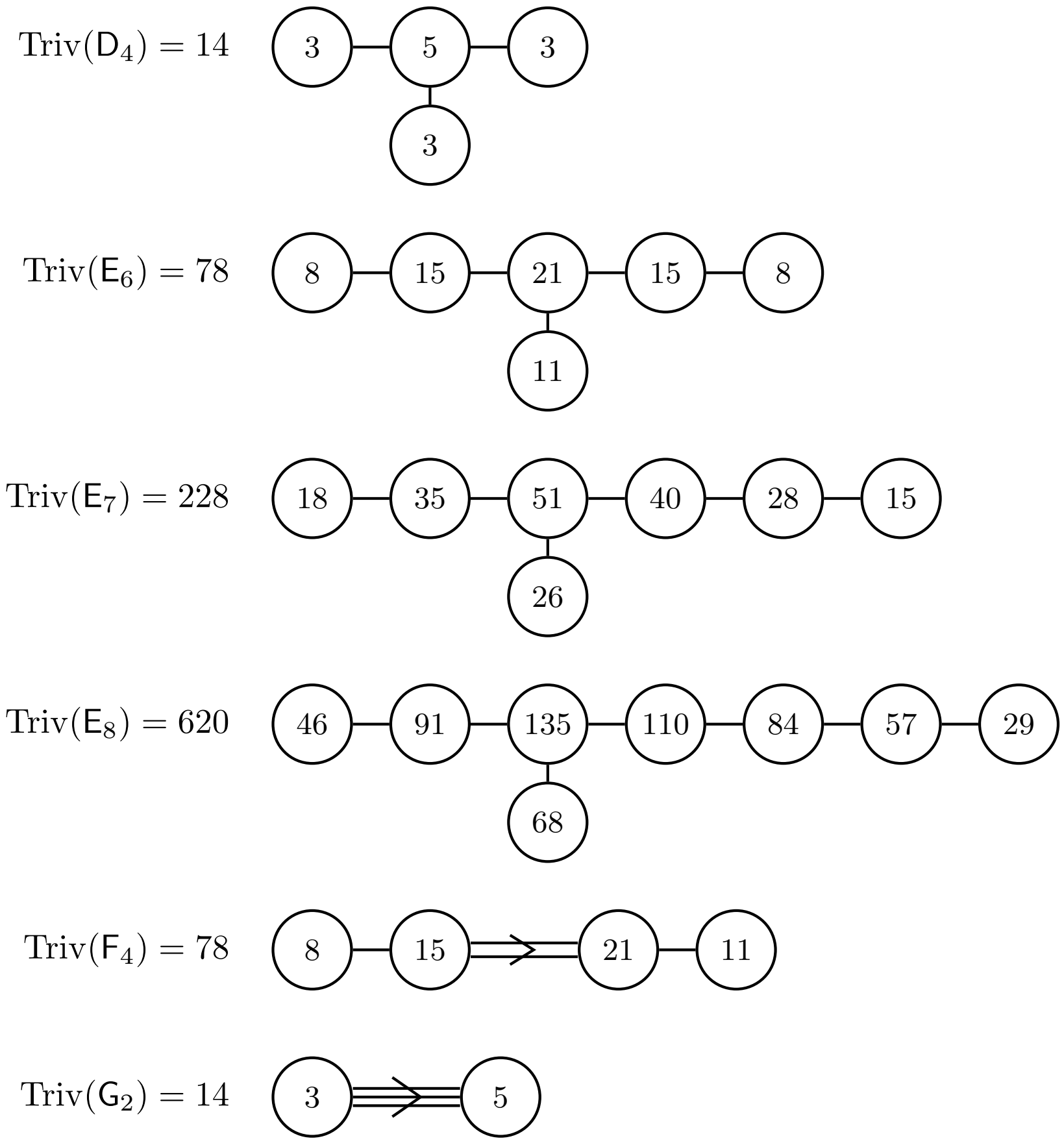
\begin{figure}[h]
\centering
\begin{tikzpicture}[x=1.2cm, y=1cm]
  \node[anchor=east] at (-0.6,0) {$\operatorname{Triv}(\mathsf{D}_4) = 14$};
  \node[dynnode] (a1) at (0,0) {$3$};
  \node[dynnode] (a2) at (1,0) {$5$};
  \node[dynnode] (a3) at (2,0) {$3$};
  \node[dynnode] (a4) at (1,-1) {$3$};
  \draw[thick] (a1) -- (a2) -- (a3);
  \draw[thick] (a2) -- (a4);
  \node[anchor=east] at (-0.6,-2.3) {$\operatorname{Triv}(\mathsf{E}_6) = 78$};
  \node[dynnode] (b1) at (0,-2.3) {$8$};
  \node[dynnode] (b2) at (1,-2.3) {$15$};
  \node[dynnode] (b3) at (2,-2.3) {$21$};
  \node[dynnode] (b4) at (3,-2.3) {$15$};
  \node[dynnode] (b5) at (4,-2.3) {$8$};
  \node[dynnode] (b6) at (2,-3.3) {$11$};
  \draw[thick] (b1) -- (b2) -- (b3) -- (b4) -- (b5);
  \draw[thick] (b3) -- (b6);
  \node[anchor=east] at (-0.6,-4.6) {$\operatorname{Triv}(\mathsf{E}_7) = 228$};
  \node[dynnode] (c1) at (0,-4.6) {$18$};
  \node[dynnode] (c2) at (1,-4.6) {$35$};
  \node[dynnode] (c3) at (2,-4.6) {$51$};
  \node[dynnode] (c4) at (3,-4.6) {$40$};
  \node[dynnode] (c5) at (4,-4.6) {$28$};
  \node[dynnode] (c6) at (5,-4.6) {$15$};
  \node[dynnode] (c7) at (2,-5.6) {$26$};
  \draw[thick] (c1) -- (c2) -- (c3) -- (c4) -- (c5) -- (c6);
  \draw[thick] (c3) -- (c7);
  \node[anchor=east] at (-0.6,-6.9) {$\operatorname{Triv}(\mathsf{E}_8) = 620$};
  \node[dynnode] (d1) at (0,-6.9) {$46$};
  \node[dynnode] (d2) at (1,-6.9) {$91$};
  \node[dynnode] (d3) at (2,-6.9) {$135$};
  \node[dynnode] (d4) at (3,-6.9) {$110$};
  \node[dynnode] (d5) at (4,-6.9) {$84$};
  \node[dynnode] (d6) at (5,-6.9) {$57$};
  \node[dynnode] (d7) at (6,-6.9) {$29$};
  \node[dynnode] (d8) at (2,-7.9) {$68$};
  \draw[thick] (d1) -- (d2) -- (d3) -- (d4) -- (d5) -- (d6) -- (d7);
  \draw[thick] (d3) -- (d8);
  \node[anchor=east] at (-0.6,-9.2) {$\operatorname{Triv}(\mathsf{F}_4) = 78$};
  \node[dynnode] (f1) at (0,-9.2) {$8$};
  \node[dynnode] (f2) at (1,-9.2) {$15$};
  \node[dynnode] (f3) at (2.6,-9.2) {$21$};
  \node[dynnode] (f4) at (3.6,-9.2) {$11$};
  \draw[thick] (f1) -- (f2);
  \draw[thick] ([yshift=2pt]f2.east) -- ([yshift=2pt]f3.west);
  \draw[thick] ([yshift=-2pt]f2.east) -- ([yshift=-2pt]f3.west);
  \draw[thick] (1.68,-9.05) -- (1.88,-9.2) -- (1.68,-9.35);
  \draw[thick] (f3) -- (f4);
  \node[anchor=east] at (-0.6,-10.7) {$\operatorname{Triv}(\mathsf{G}_2) = 14$};
  \node[dynnode] (g1) at (0,-10.7) {$3$};
  \node[dynnode] (g2) at (1.6,-10.7) {$5$};
  \draw[thick] ([yshift=2.5pt]g1.east) -- ([yshift=2.5pt]g2.west);
  \draw[thick] (g1) -- (g2);
  \draw[thick] ([yshift=-2.5pt]g1.east) -- ([yshift=-2.5pt]g2.west);
  \draw[thick] (0.68,-10.5) -- (0.92,-10.7) -- (0.68,-10.9);
\end{tikzpicture}
\caption{Vectors realizing $\operatorname{Triv}(\Phi)$ in the dominant chamber
for the types $\mathsf{D}_4$, $\mathsf{E}_6$, $\mathsf{E}_7$, $\mathsf{E}_8$, $\mathsf{F}_4$, $\mathsf{G}_2$.
Each node of the Dynkin diagram carries the coefficient of the vector on the corresponding simple root;
for $\mathsf{F}_4$ and $\mathsf{G}_2$ the arrow points from the long to the short simple roots.
The inner product of the vector with each simple root $\alpha$ is $\langle \alpha, \alpha \rangle / 2$,
except for the node of $\mathsf{E}_7$ with coefficient $15$, where it is $2$.}
\label{fig:exceptional}
\end{figure}

For $\mathsf{E}_6$, $\mathsf{E}_7$ and $\mathsf{E}_8$, we solve the integer quadratic programming problem~\eqref{eq:iqp}
computationally using the SCIP optimization suite~\cite{SCIP}.
For $\mathsf{F}_4$ and $\mathsf{G}_2$, we enumerate all lattice vectors of small norm.
The optimal vectors are given in Figure~\ref{fig:exceptional}, together with the
one for $\mathsf{D}_4$, which is the smallest case of Theorem~\ref{thm:BCD}.

The $\mathsf{E}_8$ result is confirmed by Moorhouse~\cite{Moorhouse}, who reports:
``The 310th shell is the smallest containing a regular $W(\mathsf{E}_8)$-orbit; i.e.\
the shortest lattice vectors having trivial stabilizer in $W(\mathsf{E}_8)$ have norm~$620$.''

In the realizations of the introduction, the optimal vector in the dominant chamber
is $(3, -1, -2)$ for $\mathsf{G}_2$ and $(8, 3, 2, 1)$ for $\mathsf{F}_4$,
for the simple roots $e_1 - 2e_2 + e_3$, $e_2 - e_3$ of $\mathsf{G}_2$ and
$(1,-1,-1,-1)$, $2e_4$, $e_3 - e_4$, $e_2 - e_3$ of $\mathsf{F}_4$, in the order of Figure~\ref{fig:exceptional}.

For $\mathsf{D}_4$, $\mathsf{E}_6$, $\mathsf{E}_8$, $\mathsf{F}_4$ and $\mathsf{G}_2$ these vectors are the Weyl vectors,
which are optimal and unique by Remark~\ref{rem:rho}, while for $\mathsf{E}_7$ the optimal vector is not the Weyl vector.
In all cases, the optimal vector is unique and is given by Theorem~\ref{thm:minuscule}.

\begin{remark}
We have $\operatorname{Triv}(\mathsf{G}_2) = \operatorname{Triv}(\mathsf{D}_4)$ and $\operatorname{Triv}(\mathsf{F}_4) = \operatorname{Triv}(\mathsf{E}_6)$,
and in Figure~\ref{fig:exceptional} the coefficients for $\mathsf{G}_2$ and $\mathsf{F}_4$ are those of the
nodes of $\mathsf{D}_4$ and $\mathsf{E}_6$ that are identified by the folding of their Dynkin diagrams.
This is explained in Section~\ref{sec:folding}.
\end{remark}

\section{The $\mathsf{A}_n$ case}
\label{sec:An}

The root lattice $\mathsf{A}_n$ can be represented as the set of integer vectors
$(a_1, \ldots, a_{n+1}) \in \ZZ^{n+1}$ with $\sum_{i=1}^{n+1} a_i = 0$.
The Weyl group $W(\mathsf{A}_n) \cong S_{n+1}$ acts by permuting coordinates.
A vector $v = (a_1, \ldots, a_{n+1})$ has trivial stabilizer if and only if
$a_i \neq a_j$ for all $i \neq j$.
The norm is $\lVert v \rVert^2 = \sum_{i=1}^{n+1} a_i^2$.

Thus $\operatorname{Triv}(\mathsf{A}_n)$ reduces to minimizing $\sum_{i=1}^{N} a_i^2$
over $N = n+1$ distinct integers summing to zero.
Setting $k = \lfloor N/2 \rfloor = \lfloor (n+1)/2 \rfloor$, we prove:

\begin{theorem}
\label{thm:An}
The problem $\operatorname{Triv}(\mathsf{A}_n)$ has a unique optimal solution (up to the action of $W(\mathsf{A}_n)$).
In the ambient representation, the optimal coordinates are:
\begin{enumerate}
\item For $n$ even ($N = n+1$ odd, $k = n/2$): the set $\{-k, \ldots, -1, 0, 1, \ldots, k\}$.
\item For $n$ odd ($N = n+1$ even, $k = (n+1)/2$): the set $\{-k, \ldots, -1, 1, \ldots, k\}$.
\end{enumerate}
In both cases, $\operatorname{Triv}(\mathsf{A}_n) = \dfrac{k(k+1)(2k+1)}{3}$.

Expressed in the basis of simple roots $\alpha_i = e_i - e_{i+1}$,
the optimal vector in the dominant Weyl chamber has coefficients:
\begin{enumerate}
\item For $n$ even: $c_i = i(n+1-i)/2$ for $1 \leq i \leq n$.
\item For $n$ odd: $c_i = j(n+2-j)/2$ where $j = \min(i, n+1-i)$, for $1 \leq i \leq n$.
\end{enumerate}
\end{theorem}

The proof relies on two exchange lemmas that constrain the structure of optimal solutions.

\begin{lemma}[No large gaps]
\label{lem:nogap3}
Let $S = \{a_1 < \cdots < a_N\}$ be a set of $N$ distinct integers with $\sum a_i = 0$
that minimizes $\sum a_i^2$. Then $a_{i+1} - a_i \leq 2$ for all $1 \leq i \leq N-1$.
\end{lemma}
\begin{proof}
Suppose $a_{i+1} \geq a_i + 3$ for some $i$.
Replace $a_i$ by $a_i + 1$ and $a_{i+1}$ by $a_{i+1} - 1$.
Both new values lie strictly between $a_i$ and $a_{i+1}$, so they are distinct
from each other and from all other elements of $S$.
The sum is preserved, and the change in $\sum a_i^2$ is
\begin{equation*}
(a_i+1)^2 + (a_{i+1}-1)^2 - a_i^2 - a_{i+1}^2 = 2(a_i - a_{i+1} + 1) \leq 2(-3+1) = -4 < 0,
\end{equation*}
contradicting optimality.
\end{proof}

\begin{lemma}[At most one gap]
\label{lem:onegap}
Under the same hypotheses as Lemma~\ref{lem:nogap3},
there is at most one index $i$ with $a_{i+1} - a_i = 2$.
\end{lemma}
\begin{proof}
Suppose $a_{i+1} = a_i + 2$ and $a_{j+1} = a_j + 2$ for some $i < j$.
The integers $a_i + 1$ and $a_j + 1$ are absent from $S$.
Replace $a_i$ by $a_i + 1$ and $a_{j+1}$ by $a_j + 1$.
The two new values are distinct from each other (since $a_i < a_j$ implies $a_i + 1 < a_j + 1$)
and from all remaining elements of $S$.
The sum is preserved since
$(a_i + 1) + (a_j + 1) - a_i - (a_j + 2) = 0$.
The change in $\sum a_i^2$ is
\begin{equation*}
(a_i+1)^2 + (a_j+1)^2 - a_i^2 - (a_j+2)^2 = 2(a_i - a_j - 1) < 0,
\end{equation*}
contradicting optimality.
\end{proof}

\begin{proof}[Proof of Theorem~\ref{thm:An}]
Let $S = \{a_1 < \cdots < a_N\}$ be an optimal set.
By Lemmas~\ref{lem:nogap3} and~\ref{lem:onegap}, the set $S$ is
either $N$ consecutive integers or $N+1$ consecutive integers with one element removed.

If $S$ has no gap, then $S = \{a_1, a_1+1, \ldots, a_1+N-1\}$ and
$\sum a_i = N a_1 + N(N-1)/2 = 0$, giving $a_1 = -(N-1)/2$.
This is an integer if and only if $N$ is odd, i.e., $n$ is even.
In that case $k = n/2$ and $S = \{-k, \ldots, k\}$.

If $S$ has exactly one gap, then $S = \{a_1, \ldots, a_1 + N\} \setminus \{g\}$ for some integer $g \in [a_1, a_1 + N]$.
The sum condition gives
\begin{equation*}
g = (N+1)\,a_1 + \frac{N(N+1)}{2}.
\end{equation*}
The constraint $a_1 \leq g \leq a_1 + N$ forces
$-\frac{N+1}{2} \leq a_1 \leq -\frac{N-1}{2}$.

\begin{itemize}
\item If $N = 2k$ (i.e., $n$ odd):
The only integer in that range is $a_1 = -k$, giving $g = 0$.
Thus $S = \{-k, \ldots, -1, 1, \ldots, k\}$.

\item If $N = 2k+1$ (i.e., $n$ even):
Then $a_1 \in \{-k-1, -k\}$.
For $a_1 = -k$, we get $g = k+1 = a_1 + N$, so the gap is at the top endpoint.
For $a_1 = -k-1$, we get $g = -k-1 = a_1$, so the gap is at the bottom endpoint.
In both cases, removing an endpoint reduces to the case without gap: $S = \{-k, \ldots, k\}$.
\end{itemize}

In all cases, the optimal norm is
$\sum a_i^2 = 2\sum_{j=1}^{k} j^2 = \frac{k(k+1)(2k+1)}{3}$.

The simple root coordinates $c_i = \sum_{j=1}^{i} a_j$
(where the ambient coordinates are arranged in decreasing order for the dominant chamber)
yield the stated formulas.
\end{proof}

\section{The $\mathsf{B}_n$, $\mathsf{C}_n$ and $\mathsf{D}_n$ cases}
\label{sec:BCD}

The three families are realized in $\RR^n$ with standard basis $e_1, \ldots, e_n$:
\begin{itemize}
\item $\mathsf{B}_n$ ($n \geq 2$): roots $\pm e_i$, $\pm e_i \pm e_j$; root lattice $\ZZ^n$;
Weyl group the signed permutations.
\item $\mathsf{C}_n$ ($n \geq 2$): roots $\pm 2e_i$, $\pm e_i \pm e_j$; root lattice
$\mathsf{D}_n = \{x \in \ZZ^n : \sum x_i \equiv 0 \pmod{2}\}$;
Weyl group the signed permutations.
\item $\mathsf{D}_n$ ($n \geq 4$): roots $\pm e_i \pm e_j$; root lattice $\mathsf{D}_n$;
Weyl group the signed permutations with an even number of sign changes.
\end{itemize}
We take the simple roots $\alpha_i = e_i - e_{i+1}$ for $1 \leq i \leq n-1$ together with
$\alpha_n = e_n$ for $\mathsf{B}_n$, $\alpha_n = 2e_n$ for $\mathsf{C}_n$, and $\alpha_n = e_{n-1} + e_n$ for $\mathsf{D}_n$.
The interior of the dominant chamber is
$x_1 > x_2 > \cdots > x_n > 0$ for $\mathsf{B}_n$ and $\mathsf{C}_n$, and
$x_1 > x_2 > \cdots > x_{n-1} > |x_n|$ for $\mathsf{D}_n$.
Writing $x = \sum_i c_i \alpha_i$ and $s_i = x_1 + \cdots + x_i$,
the coefficients on the simple roots are
\begin{equation}
\label{eq:coeffs}
c_i = s_i \;\text{ for all } i \;\;(\mathsf{B}_n),\qquad
c_i = s_i \;\text{ for } i \leq n-1 \text{ and } c_n = \tfrac{s_n}{2} \;\;(\mathsf{C}_n),
\end{equation}
and, for $\mathsf{D}_n$, $c_i = s_i$ for $i \leq n-2$, $c_{n-1} = (s_{n-1} - x_n)/2$ and $c_n = s_n/2$.

\begin{lemma}
\label{lem:stabBCD}
A vector $v = (x_1, \ldots, x_n)$ in the root lattice has trivial stabilizer
\begin{enumerate}
\item for $\mathsf{B}_n$ and $\mathsf{C}_n$, if and only if $|x_1|, \ldots, |x_n|$ are nonzero and pairwise distinct;
\item for $\mathsf{D}_n$, if and only if $|x_1|, \ldots, |x_n|$ are pairwise distinct.
\end{enumerate}
\end{lemma}
\begin{proof}
The stabilizer is trivial if and only if $\langle v, \alpha \rangle \neq 0$ for every root $\alpha$.
We have $\langle v, e_i \pm e_j \rangle = x_i \pm x_j$, which is nonzero for both signs
if and only if $|x_i| \neq |x_j|$,
and $\langle v, e_i \rangle$, $\langle v, 2e_i \rangle$ are nonzero if and only if $x_i \neq 0$.
\end{proof}

\begin{lemma}
\label{lem:distinct}
Let $\epsilon \in \{0, 1\}$ and let $\epsilon \leq c_1 < c_2 < \cdots < c_n$ be integers.
Then
\begin{equation*}
\sum_{i=1}^n c_i^2 \geq \sum_{i=1}^{n} (i-1+\epsilon)^2,
\end{equation*}
with equality if and only if $c_i = i - 1 + \epsilon$ for all $i$.
If $(c_1, \ldots, c_n) \neq (\epsilon, \ldots, n-1+\epsilon)$, then
\begin{equation*}
\sum_{i=1}^n c_i^2 \geq \sum_{i=1}^{n-1} (i-1+\epsilon)^2 + (n+\epsilon)^2,
\end{equation*}
with equality if and only if $c_i = i-1+\epsilon$ for $i \leq n-1$ and $c_n = n + \epsilon$.
\end{lemma}
\begin{proof}
Since the $c_i$ are strictly increasing integers, $c_i \geq c_1 + (i-1) \geq i-1+\epsilon \geq 0$,
which gives the first inequality and its equality case.
If $c_n = n-1+\epsilon$, then $c_i \leq c_n - (n-i) = i-1+\epsilon$ for all $i$,
so $c_i = i-1+\epsilon$ for all $i$.
Hence if $(c_i) \neq (\epsilon, \ldots, n-1+\epsilon)$ we have $c_n \geq n + \epsilon$,
and the second inequality follows from $c_i \geq i-1+\epsilon$ for $i \leq n-1$.
\end{proof}

Write $S(m) = \sum_{j=1}^m j^2 = m(m+1)(2m+1)/6$.

\begin{theorem}
\label{thm:BCD}
For each of $\mathsf{B}_n$, $\mathsf{C}_n$, $\mathsf{D}_n$, the problem $\operatorname{Triv}$ has a unique optimal solution
up to the action of the Weyl group.
The optimal vector $v$ in the dominant chamber, its coefficients $c_i$ on the simple roots,
and the minimum norm are the following.
\begin{enumerate}
\item $\mathsf{B}_n$, all $n \geq 2$: $v = (n, n-1, \ldots, 2, 1)$,
$c_i = i(2n+1-i)/2$ for $1 \leq i \leq n$, and
$\operatorname{Triv}(\mathsf{B}_n) = S(n)$.
\item $\mathsf{C}_n$, $n \equiv 0, 3 \pmod{4}$: $v = (n, n-1, \ldots, 2, 1)$,
$c_i = i(2n+1-i)/2$ for $i \leq n-1$, $c_n = n(n+1)/4$, and
$\operatorname{Triv}(\mathsf{C}_n) = S(n)$.
\item $\mathsf{C}_n$, $n \equiv 1, 2 \pmod{4}$: $v = (n+1, n-1, n-2, \ldots, 2, 1)$,
$c_i = i(2n+1-i)/2 + 1$ for $i \leq n-1$, $c_n = (n(n+1)+2)/4$, and
$\operatorname{Triv}(\mathsf{C}_n) = S(n) + 2n + 1$.
\item $\mathsf{D}_n$, $n \equiv 0, 1 \pmod{4}$: $v = (n-1, n-2, \ldots, 1, 0)$,
$c_i = i(2n-1-i)/2$ for $i \leq n-2$, $c_{n-1} = c_n = n(n-1)/4$, and
$\operatorname{Triv}(\mathsf{D}_n) = S(n-1)$.
\item $\mathsf{D}_n$, $n \equiv 2, 3 \pmod{4}$: $v = (n, n-2, n-3, \ldots, 1, 0)$,
$c_i = i(2n-1-i)/2 + 1$ for $i \leq n-2$, $c_{n-1} = c_n = (n(n-1)+2)/4$, and
$\operatorname{Triv}(\mathsf{D}_n) = S(n-1) + 2n - 1$.
\end{enumerate}
\end{theorem}
\begin{proof}
Let $v = (x_1, \ldots, x_n)$ have trivial stabilizer and let $c_1 < \cdots < c_n$
be the absolute values $|x_i|$ in increasing order;
they are pairwise distinct by Lemma~\ref{lem:stabBCD}.
The norm of $v$ is $\sum c_i^2$, and since changing the sign of $x_i$ changes
$\sum x_i$ by the even number $2x_i$, we have $\sum x_i \equiv \sum c_i \pmod{2}$.

For $\mathsf{B}_n$, the $c_i$ are positive, and Lemma~\ref{lem:distinct} with $\epsilon = 1$
shows that $\lVert v \rVert^2 \geq S(n)$, with equality if and only if
$\{|x_i|\} = \{1, \ldots, n\}$.
All such vectors are signed permutations of $(n, \ldots, 1)$, hence form a single orbit of $W(\mathsf{B}_n)$.

For $\mathsf{C}_n$, the $c_i$ are positive and $\sum c_i$ must be even.
If $n \equiv 0, 3 \pmod{4}$, then $\sum_{i=1}^n i = n(n+1)/2$ is even, and we conclude as for $\mathsf{B}_n$.
If $n \equiv 1, 2 \pmod{4}$, then $n(n+1)/2$ is odd, so $(c_i) \neq (1, \ldots, n)$
and the second part of Lemma~\ref{lem:distinct} gives
$\lVert v \rVert^2 \geq S(n-1) + (n+1)^2 = S(n) + 2n + 1$,
with equality if and only if $\{|x_i|\} = \{1, \ldots, n-1, n+1\}$.
This set has even sum $n(n+1)/2 + 1$, so all the corresponding vectors lie in $\mathsf{D}_n$,
and they form a single orbit of $W(\mathsf{C}_n)$.

For $\mathsf{D}_n$, at most one $c_i$ is zero, and we apply Lemma~\ref{lem:distinct} with $\epsilon = 0$.
The sum $\sum_{i=0}^{n-1} i = n(n-1)/2$ is even if and only if $n \equiv 0, 1 \pmod{4}$.
In that case the optimum is $S(n-1)$, realized exactly when $\{|x_i|\} = \{0, 1, \ldots, n-1\}$.
Otherwise the second part of Lemma~\ref{lem:distinct} gives the bound
$S(n-2) + n^2 = S(n-1) + 2n - 1$, realized exactly when $\{|x_i|\} = \{0, 1, \ldots, n-2, n\}$,
whose sum is even.
In both cases one of the $x_i$ is zero, so every signed permutation of $v$ can be
obtained with an even number of sign changes: the optimal vectors form a single orbit of $W(\mathsf{D}_n)$.

The vectors $v$ listed are the representatives of these orbits in the dominant chamber,
and the coefficients $c_i$ follow from~\eqref{eq:coeffs} and the partial sums
$\sum_{j=1}^{i} (n+1-j) = i(2n+1-i)/2$ and $\sum_{j=1}^{i} (n-j) = i(2n-1-i)/2$.
\end{proof}

\begin{remark}
The root systems $\mathsf{B}_2$ and $\mathsf{C}_2$ are the same up to a rotation and a scaling by $\sqrt{2}$,
in agreement with $\operatorname{Triv}(\mathsf{C}_2) = 10 = 2 \operatorname{Triv}(\mathsf{B}_2)$.
Also, $\operatorname{Triv}(\mathsf{C}_n) = \operatorname{Triv}(\mathsf{D}_{n+1})$ for all $n \geq 3$,
the optimal vector of $\mathsf{D}_{n+1}$ being the optimal vector of $\mathsf{C}_n$ with a coordinate $0$ appended,
and $\operatorname{Triv}(\mathsf{A}_{2n-1}) = \operatorname{Triv}(\mathsf{A}_{2n}) = 2 \operatorname{Triv}(\mathsf{B}_n)$ for all $n \geq 2$.
These equalities are explained in Section~\ref{sec:folding}.
\end{remark}

\section{Comparison with the Weyl vector}
\label{sec:weyl}

The Weyl vector of $\Phi$ is
\begin{equation*}
\rho = \frac{1}{2} \sum_{\alpha > 0} \alpha.
\end{equation*}
It lies in the interior of the dominant chamber and is characterized by
$\langle \rho, \alpha_i \rangle = \langle \alpha_i, \alpha_i \rangle / 2$
for every simple root $\alpha_i$, since the simple reflection $s_{\alpha_i}$ permutes the positive roots
other than $\alpha_i$~\cite[Proposition~1.4]{Humphreys}.
We compare it with the optimal vector $v$ in the dominant chamber given by
Theorems~\ref{thm:An} and~\ref{thm:BCD} and Figure~\ref{fig:exceptional}.

\begin{theorem}
\label{thm:weyl}
Let $\Phi$ be one of $\mathsf{A}_n$ ($n \geq 2$), $\mathsf{B}_n$ ($n \geq 2$), $\mathsf{C}_n$ ($n \geq 2$), $\mathsf{D}_n$ ($n \geq 4$),
$\mathsf{E}_6$, $\mathsf{E}_7$, $\mathsf{E}_8$, $\mathsf{F}_4$, $\mathsf{G}_2$.
Then $v = \rho$ if and only if $\rho \in L$, and otherwise $v$ is not a multiple of $\rho$.
We have $\rho \in L$ exactly for $\mathsf{A}_n$ with $n$ even, $\mathsf{C}_n$ with $n \equiv 0, 3 \pmod{4}$,
$\mathsf{D}_n$ with $n \equiv 0, 1 \pmod{4}$, $\mathsf{E}_6$, $\mathsf{E}_8$, $\mathsf{F}_4$ and $\mathsf{G}_2$.
\end{theorem}
\begin{proof}
We compare the two vectors case by case.
\begin{itemize}
\item $\mathsf{A}_n$: $\rho = \frac{1}{2}(n, n-2, \ldots, -n)$.
For $n$ even, $\rho = (k, \ldots, -k) = v$.
For $n$ odd, $\rho \notin \ZZ^{n+1}$, and $v = (k, \ldots, 1, -1, \ldots, -k)$
is not a multiple of $\rho$ since its consecutive differences are not all equal.
\item $\mathsf{B}_n$: $\rho = (n - \frac{1}{2}, \ldots, \frac{3}{2}, \frac{1}{2}) \notin \ZZ^n$,
and $v = (n, \ldots, 1)$ is not a multiple of $\rho$.
\item $\mathsf{C}_n$: $\rho = (n, \ldots, 1)$ lies in $\mathsf{D}_n$ if and only if $n(n+1)/2$ is even,
that is $n \equiv 0, 3 \pmod{4}$, and then $\rho = v$.
Otherwise $v = (n+1, n-1, \ldots, 1)$ is not a multiple of $\rho$.
\item $\mathsf{D}_n$: $\rho = (n-1, \ldots, 1, 0)$ lies in $\mathsf{D}_n$ if and only if $n(n-1)/2$ is even,
that is $n \equiv 0, 1 \pmod{4}$, and then $\rho = v$.
Otherwise $v = (n, n-2, \ldots, 1, 0)$ is not a multiple of $\rho$.
\item $\mathsf{E}_6$, $\mathsf{E}_8$, $\mathsf{F}_4$, $\mathsf{G}_2$: by Figure~\ref{fig:exceptional},
$\langle v, \alpha_i \rangle = \langle \alpha_i, \alpha_i \rangle / 2$ for every simple root,
so $v = \rho$, and $\rho \in L$.
\item $\mathsf{E}_7$: by Figure~\ref{fig:exceptional}, $\langle v, \alpha_i \rangle$ is $1$ for all
simple roots but one, where it is $2$, so $v$ is not a multiple of $\rho$.
The coefficients of $\rho$ on the simple roots are
$17, 33, 48, \frac{75}{2}, 26, \frac{27}{2}$ along the long branch and $\frac{49}{2}$
on the short one, so $\rho \notin \mathsf{E}_7$. \qedhere
\end{itemize}
\end{proof}

\begin{remark}
\label{rem:rho}
If $\rho \in L$, then $\rho$ is optimal and unique, without computation.
For $x \in L$, $\langle x, \alpha_i \rangle = \langle \rho, \alpha_i \rangle \langle x, \alpha_i^\vee \rangle$
with $\langle x, \alpha_i^\vee \rangle \in \ZZ$, so if $v \in L$ is in the open dominant chamber,
then $w = v - \rho$ is dominant. Since $\rho$ is a sum of positive roots, $\langle \rho, w \rangle \geq 0$
and $\lVert v \rVert^2 = \lVert \rho \rVert^2 + 2\langle \rho, w \rangle + \lVert w \rVert^2 \geq \lVert \rho \rVert^2$,
with equality if and only if $v = \rho$. This proves the results for $\mathsf{E}_6$, $\mathsf{E}_8$, $\mathsf{F}_4$, $\mathsf{G}_2$.
\end{remark}

The argument of Remark~\ref{rem:rho} extends to the case $\rho \notin L$,
and gives a uniform description of the optimal vector, together with its uniqueness.
Let $P = \{x \in V : \langle x, \alpha^\vee \rangle \in \ZZ \text{ for all } \alpha \in \Phi\}$
be the weight lattice~\cite[\S 2.9]{Humphreys}, with fundamental weights $\omega_1, \ldots, \omega_n$,
and let $P^+$ be the set of dominant weights.
We write $\mu \preceq \lambda$ if $\lambda - \mu$ is a nonnegative integral combination of simple roots.
A dominant weight $\lambda$ is \emph{minuscule} if it is nonzero and $\langle \lambda, \alpha^\vee \rangle \in \{0, \pm 1\}$ for all $\alpha \in \Phi$.
By Stembridge~\cite[Corollary~1.4, Proposition~1.12 and Corollary~1.13]{Stembridge},
each coset $\mu + L$ of $L$ in $P$ contains a unique dominant weight $\lambda$ with $\lambda \preceq \lambda'$
for every dominant weight $\lambda' \in \mu + L$,
and this weight is $0$ for the coset $L$ and the unique minuscule weight of the coset otherwise.

\begin{theorem}
\label{thm:minuscule}
Let $\lambda_\Phi$ be the minimal dominant weight of the coset $\rho + L$.
Then $v = \rho + \lambda_\Phi$ is the unique vector of $L$ of minimum norm in the open dominant chamber, and
\begin{equation*}
\operatorname{Triv}(\Phi) = \lVert \rho + \lambda_\Phi \rVert^2.
\end{equation*}
In particular, the optimal solution is unique up to the action of $W$.
\end{theorem}
\begin{proof}
Let $x \in L$ be in the open dominant chamber.
Since $L \subseteq P$, the integers $\langle x, \alpha_i^\vee \rangle$ are at least $1$,
and since $\langle \rho, \alpha_i^\vee \rangle = 1$, the vector $w = x - \rho$ is a dominant weight.
As $2\rho$ is the sum of the positive roots, $w \in -\rho + L = \rho + L$.
Hence $w = \lambda_\Phi + q$ with $q$ a nonnegative integral combination of simple roots, and
\begin{equation*}
\lVert x \rVert^2 = \lVert \rho + \lambda_\Phi \rVert^2 + 2 \langle \rho + \lambda_\Phi, q \rangle + \lVert q \rVert^2.
\end{equation*}
Since $\rho + \lambda_\Phi$ is in the open dominant chamber, $\langle \rho + \lambda_\Phi, q \rangle > 0$ for $q \neq 0$.
So $\lVert x \rVert^2 \geq \lVert \rho + \lambda_\Phi \rVert^2$, with equality if and only if $x = \rho + \lambda_\Phi$,
and $\rho + \lambda_\Phi$ is in $L$ and in the open dominant chamber.
\end{proof}

Remark~\ref{rem:rho} is the case $\lambda_\Phi = 0$, that is $\rho \in L$.
In the cases where $\rho \notin L$, the weight $\lambda_\Phi$ and the optimal vector of
Theorems~\ref{thm:An} and~\ref{thm:BCD} are recovered as follows:
\begin{itemize}
\item $\mathsf{A}_{2k-1}$: $\lambda_\Phi = \omega_k = \frac{1}{2}(1, \ldots, 1, -1, \ldots, -1)$
and $v = (k, \ldots, 1, -1, \ldots, -k)$;
\item $\mathsf{B}_n$: $\lambda_\Phi = \omega_n = \frac{1}{2}(1, \ldots, 1)$ and $v = (n, \ldots, 1)$;
\item $\mathsf{C}_n$, $n \equiv 1, 2 \pmod{4}$: $\lambda_\Phi = \omega_1 = e_1$ and $v = (n+1, n-1, \ldots, 1)$;
\item $\mathsf{D}_n$, $n \equiv 2, 3 \pmod{4}$: $\lambda_\Phi = \omega_1 = e_1$ and $v = (n, n-2, \ldots, 1, 0)$;
\item $\mathsf{E}_7$: $\lambda_\Phi$ is the minuscule weight $\omega_7$, dual to the simple root at the end of the long branch,
which is the node with coefficient $15$ in Figure~\ref{fig:exceptional}.
With $\lVert \rho \rVert^2 = \frac{399}{2}$, $\langle \rho, \omega_7 \rangle = \frac{27}{2}$ and $\lVert \omega_7 \rVert^2 = \frac{3}{2}$,
we get $\operatorname{Triv}(\mathsf{E}_7) = \frac{399}{2} + 27 + \frac{3}{2} = 228$.
\end{itemize}
In particular, the results for the exceptional types do not depend on the computations
of Section~\ref{sec:exceptional}, which serve as a confirmation.

\section{Folding}
\label{sec:folding}

The equalities between minimal norms observed in Sections~\ref{sec:exceptional} and~\ref{sec:BCD},
namely
\begin{gather*}
\operatorname{Triv}(\mathsf{G}_2) = \operatorname{Triv}(\mathsf{D}_4),\quad
\operatorname{Triv}(\mathsf{F}_4) = \operatorname{Triv}(\mathsf{E}_6),\quad
\operatorname{Triv}(\mathsf{C}_n) = \operatorname{Triv}(\mathsf{D}_{n+1}),\\
2\operatorname{Triv}(\mathsf{B}_n) = \operatorname{Triv}(\mathsf{A}_{2n-1}) = \operatorname{Triv}(\mathsf{A}_{2n}),
\end{gather*}
come from the automorphisms of the Dynkin diagrams of the simply laced root systems.
These diagram automorphisms and the corresponding folding of Dynkin diagrams are
at the basis of the construction of the twisted affine Lie algebras~\cite[\S 8.3]{Kac}.

Let $\Phi$ be a simply laced root system with simple roots $\Delta = \{\alpha_1, \ldots, \alpha_n\}$
and let $\sigma$ be a nontrivial automorphism of its Dynkin diagram.
Since the Gram matrix of $\Delta$ is determined by the Dynkin diagram,
$\sigma$ extends to an isometry of $V$, still denoted $\sigma$, which permutes $\Delta$.
It therefore preserves $\Phi$, the root lattice $L$ and the dominant chamber $\MC$.
Let $V^\sigma$ be the subspace of vectors fixed by $\sigma$ and $L^\sigma = L \cap V^\sigma$.
For a $\sigma$-orbit $O \subseteq \Delta$, set $\alpha_O = \sum_{\beta \in O} \beta$.
Since $\sigma$ permutes the basis $\Delta$ of $L$, the vectors $\alpha_O$ form a basis of $L^\sigma$.
Two roots in the same orbit are orthogonal, except for the middle orbit $\{\alpha_n, \alpha_{n+1}\}$
of $\mathsf{A}_{2n}$, where they are adjacent;
so $\langle \alpha_O, \alpha_O \rangle = 2|O|$, except for that orbit, where it is $2$.
A direct computation of the Gram matrix of the $\alpha_O$ shows that they form a system of simple roots
of a root system $\Phi^\sigma$ in $V^\sigma$, whose type and root lattice $L^\sigma = \ZZ\Phi^\sigma$ are given in the following table.
\begin{center}
\begin{tabular}{ccccc}
$\Phi$ & order of $\sigma$ & $\Phi^\sigma$ & $\lVert \alpha_O \rVert^2$ & $L^\sigma$\\
\hline
$\mathsf{A}_{2n-1}$, $n \geq 2$ & $2$ & $\sqrt{2}\,\mathsf{B}_n$ & $2, 4$ & $\sqrt{2}\,\ZZ^n$\\
$\mathsf{A}_{2n}$, $n \geq 2$ & $2$ & $\sqrt{2}\,\mathsf{B}_n$ & $2, 4$ & $\sqrt{2}\,\ZZ^n$\\
$\mathsf{D}_{n+1}$, $n \geq 3$ & $2$ & $\mathsf{C}_n$ & $2, 4$ & $\mathsf{D}_n$\\
$\mathsf{D}_4$ & $3$ & $\mathsf{G}_2$ & $2, 6$ & $\mathsf{A}_2$\\
$\mathsf{E}_6$ & $2$ & $\mathsf{F}_4$ & $2, 4$ & $\mathsf{D}_4$
\end{tabular}
\end{center}
Since a root system is determined up to isometry by the Gram matrix of its simple roots,
$\Phi^\sigma$ is isometric to the realization of the introduction,
except for $\mathsf{B}_n$, which appears scaled by $\sqrt{2}$, so that $\operatorname{Triv}(\sqrt{2}\,\mathsf{B}_n) = 2\operatorname{Triv}(\mathsf{B}_n)$.
For $\mathsf{A}_{2n-1}$ and $\mathsf{D}_{n+1}$, the subalgebra of the simple Lie algebra of type $\Phi$
fixed by the automorphism induced by $\sigma$ has type $\mathsf{C}_n$ and $\mathsf{B}_n$ respectively~\cite[Proposition~8.3~e)]{Kac},
the dual of $\Phi^\sigma$.
A root system and its dual have the same reflection hyperplanes, hence the same regular vectors;
it is the lattice $L^\sigma$ that selects $\Phi^\sigma$ among the two.

\begin{proposition}
\label{prop:folding}
Let $\Phi$ and $\sigma$ be as above.
Then $\operatorname{Triv}(\Phi^\sigma) \geq \operatorname{Triv}(\Phi)$,
with equality if and only if some vector of $L$ in the open dominant chamber of norm $\operatorname{Triv}(\Phi)$ is fixed by $\sigma$.
If the optimal solution for $\Phi$ is unique up to the action of $W(\Phi)$, then equality holds,
the optimal solution for $\Phi^\sigma$ is unique up to the action of $W(\Phi^\sigma)$,
and the optimal vector $v$ in the dominant chamber is the same for $\Phi$ and $\Phi^\sigma$.
Its coefficient on $\alpha_O$ is the common coefficient of $v$ on the roots of $O$.
\end{proposition}
\begin{proof}
For $x \in V^\sigma$ and $\beta \in O$, we have $\langle x, \beta \rangle = \langle x, \sigma(\beta) \rangle$,
so $\langle x, \beta \rangle$ is constant on $O$ and $\langle x, \alpha_O \rangle = |O| \langle x, \beta \rangle$.
Hence the open dominant chamber of $\Phi^\sigma$, defined by $\langle x, \alpha_O \rangle > 0$ for all orbits $O$,
is $\MC^\circ \cap V^\sigma$, where $\MC^\circ$ is the open dominant chamber of $\Phi$.
Since $\operatorname{Triv}$ can be computed in the open dominant chamber,
\begin{equation*}
\operatorname{Triv}(\Phi^\sigma) = \min \left\{ \lVert x \rVert^2 : x \in L^\sigma \cap \MC^\circ \right\}
\geq \min \left\{ \lVert x \rVert^2 : x \in L \cap \MC^\circ \right\} = \operatorname{Triv}(\Phi),
\end{equation*}
with equality if and only if a minimizer of the right-hand side lies in $V^\sigma$.

Every $W(\Phi)$-orbit meets $\MC^\circ$ in at most one point,
so if the optimal solution for $\Phi$ is unique up to $W(\Phi)$, there is a unique minimizer $v \in L \cap \MC^\circ$.
Since $\sigma$ preserves $L$, $\MC^\circ$ and the norm, $\sigma(v)$ is also a minimizer, so $\sigma(v) = v$ and equality holds.
Any minimizer for $\Phi^\sigma$ in $L^\sigma \cap \MC^\circ$ is then a minimizer for $\Phi$, hence equal to $v$.
Finally, writing $v = \sum_i c_i \alpha_i$, the relation $\sigma(v) = v$ shows that $c_i$ is constant
on each orbit $O$, with value $c_O$ say, and $v = \sum_O c_O \alpha_O$.
\end{proof}

By Theorem~\ref{thm:minuscule}, or by Theorems~\ref{thm:An} and~\ref{thm:BCD} and Remark~\ref{rem:rho},
the optimal solutions for $\mathsf{A}_n$, $\mathsf{D}_n$ and $\mathsf{E}_6$ are unique,
so Proposition~\ref{prop:folding} gives the equalities above.
It also gives a second proof of the values of $\operatorname{Triv}$ for $\mathsf{B}_n$, $\mathsf{C}_n$, $\mathsf{F}_4$ and $\mathsf{G}_2$,
and of the uniqueness of their optimal solutions.
The last assertion of Proposition~\ref{prop:folding} explains the coefficients in Figure~\ref{fig:exceptional},
and the fact that the coefficients $c_i = i(2n+1-i)/2$ of the optimal vector of $\mathsf{B}_n$
are those of $\mathsf{A}_{2n-1}$ and $\mathsf{A}_{2n}$ for $i \leq n$ in Theorem~\ref{thm:An}.
For $\mathsf{D}_{n+1}$, the space $V^\sigma$ is the hyperplane $x_{n+1} = 0$,
which explains why the optimal vector of $\mathsf{D}_{n+1}$ is that of $\mathsf{C}_n$ with a coordinate $0$ appended.
Note that $\mathsf{D}_4$ admits diagram automorphisms of order $2$ and $3$, giving
$\operatorname{Triv}(\mathsf{C}_3) = \operatorname{Triv}(\mathsf{G}_2) = \operatorname{Triv}(\mathsf{D}_4) = 14$,
while $\mathsf{E}_7$ and $\mathsf{E}_8$ have no nontrivial diagram automorphism.

\section*{Acknowledgments}

The author thanks Klaus Hulek for suggesting the problem.

\section*{Tool and computational resource disclosure}

Following the Leiden Declaration on Artificial Intelligence and Mathematics~\cite{Leiden},
we disclose the automated tools used in this work.
The optimal vectors for the exceptional types were computed with programs written in GAP,
using the SCIP optimization suite~\cite{SCIP} for $\mathsf{E}_6$, $\mathsf{E}_7$ and $\mathsf{E}_8$
and the enumeration of short lattice vectors for $\mathsf{F}_4$ and $\mathsf{G}_2$.
The AI coding assistant Claude Code (Anthropic), with the large language model Claude Opus~5.5,
was used to draft and edit parts of the text and of the proofs, including the treatment of the types $\mathsf{B}_n$, $\mathsf{C}_n$,
$\mathsf{F}_4$ and $\mathsf{G}_2$, the comparison with the Weyl vector, Theorem~\ref{thm:minuscule} and Section~\ref{sec:folding},
to search the literature, and to check the results numerically by independent enumerations.
All statements and proofs were checked by the author, who takes full responsibility for the content of the paper.

\end{document}